\documentclass[11pt]{article}
\usepackage{amsfonts}
\usepackage{color}
\newtheorem{theo}{Theorem}[section]
\newtheorem{lem}[theo]{Lemma}
\newtheorem{cor}[theo]{Corollary}

\newcommand{\mysection}[1]{\section{#1} \setcounter{equation}{0}}
\newcommand{\proof}{{\sc Proof.} \quad}
\newcommand{\proofc}{{\sc Proof} \ }
\newcommand{\be}{\begin{equation} \label}
\newcommand{\ee}{\end{equation}}
\newcommand{\bea}{\begin{eqnarray}\label}
\newcommand{\eea}{\end{eqnarray}}
\newcommand{\bas}{\begin{eqnarray*}}
\newcommand{\eas}{\end{eqnarray*}}
\newcommand{\bit}{\begin{itemize}}
\newcommand{\eit}{\end{itemize}}
\newcommand{\qed}{\hfill$\Box$ \vskip.2cm}
\newcommand{\nn}{\nonumber}
\newcommand{\R}{\mathbb{R}}

\newcommand{\pO}{\partial\Omega}

\newcommand{\abs}{\\[5pt]}

\newcommand{\tm}{T_{max}}
\newcommand{\io}{\int_\Omega}

\begin{document}
\title{ Boundedness in a two-dimensional chemotaxis-haptotaxis system}
\author{
Youshan Tao\footnote{taoys@dhu.edu.cn}\\
{\small Department of Applied Mathematics, Dong Hua University,}\\
{\small Shanghai 200051, P.R.~China} }
\date{}
\maketitle
\begin{abstract}
\noindent This work studies the chemotaxis-haptotaxis system
 \bas
 \left\{ \begin{array}{ll}
    u_t= \Delta u - \chi \nabla \cdot (u\nabla v) - \xi \nabla \cdot (u\nabla w) + \mu u(1-u-w),
    &\qquad x\in \Omega, \, t>0, \\[1mm]
    v_t=\Delta v-v+u,
    &\qquad x\in \Omega, \, t>0, \\[1mm]
    w_t=-vw,
    &\qquad x\in \Omega, \, t>0,
    \end{array} \right.
 \eas
in a bounded smooth domain $\Omega\subset\R^2$ with zero-flux
boundary conditions, where the parameters $\chi, \xi$ and $\mu$ are
assumed to be positive. It is shown that under appropriate
regularity assumption on the initial data $(u_0, v_0, w_0)$, the
corresponding initial-boundary problem possesses a unique classical
solution which is global in time and bounded. In addition to coupled
estimate techniques, a novel ingredient in the proof is to establish
a one-sided pointwise estimate, which connects $\Delta w$ to $v$ and
thereby enables us to derive useful energy-type inequalities that
bypass $w$. However, we note that the approach developed in this
paper seems to be confined to the two-dimensional setting.\abs
\end{abstract}
\newpage
\mysection{Introduction}
\subsection{Chemotaxis-haptotaxis model}
We consider the chemotaxis-haptotaxis system
 \be{0}
    \left\{ \begin{array}{ll}
    u_t= \Delta u - \chi \nabla \cdot (u\nabla v) - \xi \nabla \cdot (u\nabla w) + \mu u(1-u-w),
    &\qquad x\in \Omega, \, t>0, \\[1mm]
    v_t=\Delta v-v+u,
    &\qquad x\in \Omega, \, t>0, \\[1mm]
    w_t=-vw,
    &\qquad x\in \Omega, \, t>0,
    \end{array} \right.
\ee in a physical smoothly bounded domain $\Omega\subset \R^n$,
$n\in \{2, 3\}$, under zero-flux boundary conditions
 \be{01}
     \frac{\partial u}{\partial\nu}-\chi u\frac{\partial v}{\partial\nu}-\xi u \frac{\partial
    w}{\partial\nu}=\frac{\partial v}{\partial\nu}=0,
    \qquad x\in \partial\Omega, \, t>0
 \ee
and with prescribed initial data
 \be{02}
          u(x,0)=u_0(x), \quad v(x, 0)=v_0(x), \quad w(x,0)=w_0(x),
    \qquad x\in \Omega,
 \ee
where $\frac{\partial}{\partial \nu}$ denotes differentiation with
respect to outward normal on $\pO$,  and the parameters $\chi, \xi$
and $\mu$ are assumed to be positive. This system was initially
proposed by Chaplain and Lolas \cite{chaplain_lolas_2005,
chaplain_lolas_2006} to model the process of cancer cell invasion of
surrounding tissue. In this context, $u$ represents the density of
cancer cell, $v$ denotes the concentration of enzyme, and $w$ stands
for the density of extracellular matrix (tissue). In addition to
random motion, cancer cells bias their movement towards a gradient
of diffusible enzyme as well as a gradient of non-diffusible tissue
by detecting matrix molecules such as vitronectin adhered therein.
We refer to the aforementioned two directed migrations of cancer
cells as chemotaxis and haptotaxis, respectively. The cancer cells
are also assumed to undergo birth and death in a logistic manner,
competing for space with healthy tissue. The enzyme is produced by
cancer cells, and it is supposed to be influenced by diffusion and
degradation. The tissue is stiff in the sense that it does not
diffuse, but it could be degraded by enzyme upon contact.
\subsection{Previous related works on global well-posedness}
When $w\equiv 0$, the PDE system (\ref{0}) is reduced to the {\it
chemotaxis-only system}
 \bas
    \left\{ \begin{array}{ll}
    u_t= \Delta u - \chi \nabla \cdot (u\nabla v) + \mu u(1-u),
    &\qquad x\in \Omega, \, t>0, \\[1mm]
    v_t=\Delta v-v+u,
    &\qquad x\in \Omega, \, t>0.
    \end{array} \right.
\eas This system has been widely studied. In the case $\mu=0$,
solutions may blow up in finite time when $n\ge 2$
\cite{herrero_velazquez1997_pisa, nagai2001, win_arX2011}; however,
it is known that arbitrarily small $\mu>0$ guarantee the global
existence and boundedness of solutions when $n=2$ \cite{OTYM}, and
that appropriately large $\frac{\mu}{\chi}$ preclude blow-up in the
case $n\ge 3$ \cite{win_cpde2010}. Very recently, it is asserted
that sufficiently large $\frac{\mu}{\chi}$ enforce the stability of
constant equilibria \cite{win_jde2014}.\abs
When $\chi=0$,  the PDE system (\ref{0}) becomes the {\it
haptotaxis-only system}
 \bas
    \left\{ \begin{array}{ll}
    u_t= \Delta u - \xi \nabla \cdot (u\nabla w) + \mu u(1-u-w),
    &\qquad x\in \Omega, \, t>0, \\[1mm]
    v_t=\Delta v-v+u,
    &\qquad x\in \Omega, \, t>0, \\[1mm]
    w_t=-vw,
    &\qquad x\in \Omega, \, t>0.
    \end{array} \right.
\eas Global existence theories for this system were explored in
\cite{Corrias1, Corrias2, WW}, whereas the boundedness and
asymptotic behavior of solution was studied in \cite{LM2}. The above
results rule out the possibility of blow-up of solutions to this
haptotaxis-only system, although the solutions may exhibit some
pattern for certain ranges of model parameters and initial data
\cite{chaplain_lolas_2006}.\abs
Compared with the chemotaxis-only system and the haptotaxis-only
system, the coupled chemotaxis-haptotaxis system (\ref{0}) is much
less understood. As far as we know, when $n\in \{2, 3\}$, the global
boundedness of solutions for this system has remained pending so
far, although the global existence was examined in \cite{tao_jmaa09,
tao_wang_non}. The purpose of this work is to answer this issue of
boundedness when $n=2$.
\subsection{Main results}
As to the above initial data we suppose that for some $\alpha \in
(0,1)$ we have \be{init}
    \left\{ \begin{array}{l}
    u_0 \in C^{0}(\overline{\Omega}) \quad \mbox{with} \quad u_0 \ge 0 \quad \mbox{in } \Omega
    \quad \mbox{and} \quad u_0 \not\equiv 0, \\[1mm]
    v_0 \in W^{1, \infty}(\Omega) \quad \mbox{with} \quad v_0 \ge 0 \quad \mbox{in }
    \Omega, \\[1mm]
    w_0 \in C^{2+\alpha}(\overline{\Omega}) \quad \mbox{with} \quad  w_0 > 0 \quad \mbox{in } \overline{\Omega}
        \quad \mbox{with} \quad \frac{\partial w_0}{\partial\nu}=0 \quad \mbox{on } \pO.
    \end{array} \right.
\ee
Under these assumptions, our main result reads as follows.
\begin{theo}\label{theo1}
  Let $n=2$, and suppose that $\chi>0$, $\xi>0$ and $\mu>0$. Then for each $(u_0, v_0, w_0)$ fulfilling (\ref{init}),
 (\ref{0})-(\ref{02}) possesses a unique classical solution  which is global in time and bounded in $\Omega \times
  (0,\infty)$.
\end{theo}
To the best of our knowledge, this is the first boundedness result
addressing the full parabolic-parabolic-ODE chemotaxis-haptotaxis
model, despite there exist some boundedness and stabilization
results on the {\it simplified} parabolic-elliptic-ODE
chemotaxis-haptotaxis model in which the spatiotemporal evolution of
chemical concentration is described by the elliptic equation
$0=\Delta v-v+u$ replacing the original parabolic counterpart
(\cite{taowin_non, taowin8b}). Here we should point out that the
full chemotaxis-haptotaxis system is much more mathematically
challenging than the aforementioned simplified one.

\subsection{Approaches used in the paper}
A main technical difficulty in the proof of Theorem \ref{theo1}
emanates from consequences of the strong coupling in (\ref{0}) on
the spatial regularity of $u$, $v$ and $w$. Another analytical
obstacle stems from the fact that a bound of $\nabla v$ in
$L^p(\Omega)$ cannot leads to a time-independent bound for $\nabla
w$ in $L^p(\Omega)$, because
 \bas
\nabla w(x,t) = \nabla w_0(x) \, e^{-\int_0^t v(x,s)ds}
    - w_0(x) \, e^{-\int_0^t v(x,s)ds} \, \int_0^t \nabla v(x,s) ds
 \eas
in which the last term $\int_0^t \nabla v(x,s) ds$ is non-local in
time. It will be crucial to our approach to build a one-sided
pointwise estimate which connects $\Delta w$ to $v$ (see Lemma
\ref{lem_deltaw} below). Relying on such a pointwise estimate, we
can derive two useful energy-type inequalities that bypass $w$ (see
Lemmata \ref{lem2.3} and \ref{lem3.3} below). Using such information
along with coupled estimate techniques, we establish estimates on
$\io u^2+\io |\nabla v|^4$ and $\io u^p$ for any $p>2$, which
results in the boundedness of $u$ in $L^\infty(\Omega)$ by
performing the Moser iteration procedure (see Lemmata \ref{lem3.8},
\ref{lem3.10} and \ref{lem4.1} below).\abs
Finally, we mention that the methods developed in this paper are
restricted to the two-dimensional setting, and thus the boundedness
for the corresponding three-dimensional problem largely remains
open.

\mysection{Local existence and a one-sided pointwise estimate for
$\Delta w$}
With a slight adaption to the proof of \cite[Lemma 4.1]{taowin8b},
we have the following statement on local existence.
\begin{lem}\label{lem_loc}
  Let $\chi> 0, ~\xi> 0$ and $\mu > 0$. Then for any $u_0, v_0$ and $w_0$ fulfilling (\ref{init})
  there exists $\tm \in (0,\infty]$ with the property that (\ref{0})-(\ref{02}) possesses a unique classical solution
  \bas
    && u \in C^0(\bar\Omega \times [0,\tm)) \cap C^{2,1}(\bar\Omega \times (0,\tm)), \nn\\
    && v \in C^0(\bar\Omega \times [0,\tm)) \cap C^{2,1}(\bar\Omega \times (0,\tm)), \nn\\
    && w \in C^{2,1}(\bar\Omega \times [0,\tm)),
  \eas
  such that
  \be{bound}
    u\ge 0,\quad v\ge 0 \quad \mbox{and} \quad 0< w\le \|w_0\|_{L^\infty(\Omega)}
    \qquad \mbox{for all } (x,t) \in \Omega\times [0,\tm)
  \ee
  and such that
  \be{ext_crit}
        \mbox{either $\tm=\infty$, \qquad or} \qquad
    \|u(\cdot,t)\|_{L^\infty(\Omega)} \to \infty
        \quad \mbox{as  $t\nearrow \tm$}.
  \ee
\end{lem}
Since the third equation in (\ref{0}) is an ODE, $w$ can be
expressed explicitly in terms of $v$. This results in the
representation formulae
\begin{eqnarray}
    w(x,t) &=& w_0(x) \, e^{-\int_0^t v(x,s)ds} \qquad \mbox{and}  \label{w}\\
    \nabla w(x,t) &=& \nabla w_0(x) \, e^{-\int_0^t v(x,s)ds}
    - w_0(x) \, e^{-\int_0^t v(x,s)ds} \, \int_0^t \nabla v(x,s) ds \label{nablaw}
\eea as well as \bea{deltaw}
    -\Delta w(x,t)
        &=& -\Delta w_0(x) \, e^{-\int_0^t v(x,s)ds}
        +2 e^{-\int_0^t v(x,s)ds} \nabla w_0(x) \cdot \int_0^t \nabla v(x,s)ds \nn\\
    & & - w_0(x) \, e^{-\int_0^t v(x,s)ds} \cdot \bigg| \int_0^t \nabla v(x,s)ds \bigg|^2
        + w_0(x) \, e^{-\int_0^t v(x,s)ds} \int_0^t \Delta v(x,s)ds
\eea for $(x,t) \in \Omega \times (0,\tm)$.\\
The following one-sided pointwise estimate for $-\Delta w$ will
serve as a cornerstone for our subsequent analysis (see the proofs
of Lemmata \ref{lem2.3} and \ref{lem3.3} below) .
\begin{lem}\label{lem_deltaw}
 Assume that $\chi>0, \xi>0$ and $\mu>0$, and let $(u,v,w)$ solve (\ref{0})-(\ref{02}) in $\Omega \times (0,T)$
 with some some $(u_0, v_0, w_0)$
  satisfying (\ref{init}). Then
  \bea{w1}
    -\Delta w(x,t)
    \le \|w_0\|_{L^\infty(\Omega)} \cdot v(x, t) + K \qquad
    \textrm{ for all $x\in\Omega$ and $t\in (0,T)$},
  \eea
  where
  \be{K}
K :=\|\Delta w_0\|_{L^\infty(\Omega)} +4\|\nabla
\sqrt{w_0}\|_{L^\infty(\Omega)}^2
+\frac{\|w_0\|_{L^\infty(\Omega)}}{e}.
  \ee
\end{lem}
\proof Start from (\ref{deltaw}), the nonnegativity of $v$ leads to
 \be{28}
-\Delta w_0(x) \, e^{-\int_0^t v(x,s)ds}\le \|\Delta
w_0\|_{L^\infty(\Omega)} \qquad \textrm{ for all $x\in\Omega$ and
$t\in (0,T)$}
 \ee
and a simple but important observation yields
 \bea{29}
 e^{-\int_0^t v(x,s)ds} &\bigg[ &2\nabla w_0(x) \cdot \int_0^t \nabla v(x,s)ds
    - w_0(x)  \cdot \bigg| \int_0^t \nabla v(x,s)ds
    \bigg|^2\bigg]\nn\\[1mm]
   &=&  -w_0(x)e^{-\int_0^t v(x,s)ds} \, \bigg|\int_0^t \nabla v(x,s)ds -\frac{\nabla w_0(x)}{w_0(x)}\bigg|^2
        +e^{-\int_0^t v(x,s)ds} \, \frac{|\nabla
        w_0(x)|^2}{w_0(x)}\nn\\
   &\le& e^{-\int_0^t v(x,s)ds} \, \frac{|\nabla
        w_0(x)|^2}{w_0(x)}\nn\\
  &\le& \frac{|\nabla   w_0(x)|^2}{w_0(x)}\nn\\
  &\le& 4\|\nabla
\sqrt{w_0}\|_{L^\infty(\Omega)}^2 \qquad \textrm{ for all
$x\in\Omega$ and $t\in (0,T)$}
 \eea
thanks to the nonnegativity of $v$ and the positivity of $w_0$.

We now turn to estimate the last term in (\ref{deltaw}).  By the
second equation in (\ref{0}) and the nonnegativity of $u, v, v_0$
and $w_0$ we have
  \bea{210}
    w_0(x) \, e^{-\int_0^t v(x,s)ds} \int_0^t \Delta v(x,s)ds
    &=& w_0(x) \, e^{-\int_0^t v(x,s)ds} \, \int_0^t \Big(v_t(x, s) +v(x,s)-u(x,s)\Big)ds
    \nn\\[1mm]
    &\le& w_0(x) \, e^{-\int_0^t v(x,s)ds} \, \bigg(v(x, t)-v_0(x)+\int_0^t v(x,s)ds
    \bigg)\nn\\[1mm]
    &\le& \|w_0\|_{L^\infty(\Omega)} \cdot v(x, t)  +\frac{\|w_0\|_{L^\infty(\Omega)}}{e}
      \eea
 for all $x\in\Omega$ and  $t\in (0,T)$, where we have used the facts that $ze^{-z} \le \frac{1}{e}$ for all $z\in\R$
 and that $0<e^{-\int_0^t v(x,s)ds}\le 1$ thanks to $v\ge 0$.  Finally,
 collecting (\ref{28})-(\ref{210}) in conjunction with (\ref{deltaw}) yields (\ref{w1}).
\qed
Here we stress the fact that the pointwise estimate (\ref{w1})
connects $\Delta w$ to $v$, which enables us to establish the
following useful energy-type inequality that will be used in the
proofs of Lemmata \ref{lem3.5}, \ref{lem3.10} and \ref{lem4.1}
below.

\begin{lem}\label{lem2.3}
Let $n\in \{2, 3\}$, $T\in (0, \tm)$, $\chi>0,~ \xi>0$ and $\mu>0$ ,
and assume (\ref{init}). Then the solution of (\ref{0})-(\ref{02})
satisfies
\bea{2.11}
      \frac{1}{p}\frac{d}{dt} \io u^p +\frac{p-1}{2}\io  u^{p-2}|\nabla u|^2
      & \le &\frac{(p-1) \chi^2}{2}
       \io  u^p|\nabla v|^2 +\xi \|w_0\|_{L^\infty(\Omega)} \io u^p v \nn\\
       & &+(\mu+\xi K)\io  u^p -\mu \io u^{p+1}
       \eea
       for any $p>1$ and each $t\in (0, T)$.
\end{lem}
\proof We test the first equation in (\ref{0}) by $u^{p-1}$, which
leads to the identity
 \bea{2.12}
\frac{1}{p} \frac{d}{dt} \io u^p +(p-1)\io u^{p-2}|\nabla u|^2
 &=&\chi (p-1)\io u^{p-1}\nabla u\cdot \nabla v+\xi (p-1)\io u^{p-1}\nabla
u\cdot \nabla w \nn\\
& & +\mu\io u^p -\mu\io u^{p+1} -\mu \io u^pw
 \eea
for all $t\in (0, T)$. Here using the Young inequality we see that
 \be{2.13}
\chi(p-1) \io u\nabla u\cdot \nabla v \le \frac{p-1}{2} \io u^{p-2}
|\nabla u|^2+\frac{(p-1)\chi^2}{2}\io u^2 \nabla v|^2 \qquad
\textrm{for all $t\in (0, T)$}.
 \ee
Unlike the handling of the above chemotaxis-related integral, we now
invoke Lemma \ref{lem_deltaw} in conjunction with integration by
parts and the Young inequality to estimate the haptotaxis-related
integral on the right of (\ref{2.12})
 \bea{2.14}
\xi (p-1)\io u^{p-1}\nabla u\cdot \nabla w
 &=& -\frac{(p-1)\xi}{p} \io u^p \Delta w\nn\\
 &\le& \frac{(p-1)\xi}{p}\, \|w_0\|_{L^\infty(\Omega)}  \io u^p v + \frac{(p-1)\xi
 K}{p}\io u^p\nn\\
 &\le& \xi \|w_0\|_{L^\infty(\Omega)} \io u^p v + \xi K \io u^p
 \qquad \textrm{for all
$t\in (0, T)$}
 \eea
thanks to $0<\frac{p-1}{p}<1$. We next observe that
 \be{2.15}
 -\mu \io u w \le 0 \qquad
\textrm{for all $t\in (0, T)$},
 \ee
because $w\ge 0$. Finally, collecting (\ref{2.12})-(\ref{2.15})
yields (\ref{2.11}) upon a simple rearrangement.\qed
The following basic property on mass can be easily checked.
\begin{lem}\label{mass}
The solution $(u, v, w)$ of (\ref{0})-(\ref{02}) fulfills
 \begin{eqnarray}
      \io u(x, t) dx &\le& m^*:=\max\Big\{ |\Omega|, \io u_0(x) dx\Big\} \qquad \mbox{for all } t \in (0,\tm). \label{2.16}
  \eea
\end{lem}
\proof
  We integrate the first equation in (\ref{0}) with respect to space to obtain
  \bea{2.17}
    \frac{d}{dt} \io u
    = \mu \io u - \mu \io u^2 - \mu \io uw
    \le \mu \io u - \mu \io u^2
    \qquad \mbox{for all } t\in (0,\tm),
  \eea
  because $w> 0$ by Lemma \ref{lem_loc}.
  Thanks to the Cauchy-Schwarz inequality, we have $\io u^2 \ge \frac{1}{|\Omega|} (\io u)^2$, and thus
  (\ref{2.17}) entails that $y(t):=\io u(x,t)dx, \ t\in [0,\tm)$, satisfies
  \bas
    y'(t) \le \mu y(t) - \frac{\mu}{|\Omega|} y^2(t)
    \qquad \mbox{for all } t\in (0,\tm).
  \eas
  By an ODE comparison, we therefore obtain that $y(t)\le \max \{|\Omega|, y(0)\}$, which precisely results in (\ref{2.16}).
  \qed
\mysection{Energy-type estimates for any $\mu>0$}
In order to establish a bound of $u$ in $L^\infty(\Omega)$, we first
build an estimate on $\io u\ln u$ as a starting point of our
reasoning, which heavily depends on the pointwise estimate
(\ref{w1}) for $-\Delta w$.
\subsection{A coupled estimate on $\io u\ln u +\io |\nabla v|^2$}
We begin with an elementary lemma.
\begin{lem}\label{lem3.1}

  Let $\mu>0$ and $A>0$. Then there exists $L:= L(\mu, A)>0$ such that
  \be{3.1}
   (1+ \mu) z \ln z +A z^2- \mu z^2 \ln z \le L
    \qquad \mbox{for all } z>0.
  \ee
\end{lem}
\proof
  The function $\varphi: \ [0,\infty) \to \R$ defined by
  \bas
    \varphi(z):=\left\{ \begin{array}{ll}
    (1+\mu) z\ln z +A z^2- \mu z^2 \ln z, \qquad & z>0, \\[1mm]
    0, & z=0,
    \end{array} \right.
  \eas
  satisfies
  \bas
    \frac{\varphi(z)}{z^2 \ln z} \to -\mu \qquad \mbox{as } z\to\infty,
  \eas
  so that for some $z_0>0$ we have $\varphi<0$ on $(z_0,\infty)$. Since clearly $\varphi$ is continuous
  on $[0,\infty)$, (\ref{3.1}) thus holds with $L:=\max_{z\in [0,z_0]} \varphi(z)$.
\qed
In the two-dimensional setting, using the properties of the Neumann
heat semigroup (\cite{win_jde}) and the estimate (\ref{2.16}) on $u$
in $L^1(\Omega)$ provided by Lemma \ref{mass}, we can derive a $L^p$
estimate on $v$.
\begin{lem}\label{lem3.2}
Let $n=2$, and assume (\ref{init}). Then for any $1\le p<\infty$
there exists a positive constant $M(p):=M(p, |\Omega|,
\|u_0\|_{L^1(\Omega)}, \|v_0\|_{L^\infty(\Omega})>0$ such that the
solution of (\ref{0})-(\ref{02}) satisfies
 \be{3.2}
\io v^p \le M(p) \qquad \mbox{for all $t\in (0, \tm)$}.
 \ee
\end{lem}
\proof Since the proof was given in \cite[Lemma 3.1]{liutao}, we
refrain us from repeating it here. \qed
Strongly depending on the estimate (\ref{w1}) for $-\Delta w$ once
again, along with the above two preparations, we now can establish
two estimates on $\io u\ln u$ and $\io |\nabla v|^2$ via coupled
estimate techniques.
\begin{lem}\label{lem3.3}
Let $n=2$, $T\in (0, \tm)$, $\chi>0,~ \xi>0$ and $\mu>0$ , and
assume (\ref{init}). Then there exists $C>0$ independent of $T$ such
that the solution of (\ref{0})-(\ref{02}) possesses the properties
\begin{eqnarray}
      \io u\ln u &\le& C \qquad \mbox{for all } t \in (0,T)\quad \mbox{and}\label{3.3}
      \\[1mm]
      \io |\nabla v|^2 &\le& C \qquad \mbox{for all } t \in
      (0,T). \label{3.4}
  \eea
\end{lem}
\proof Testing the first equation of (\ref{0}) against $(1+\ln u)$,
we obtain
 \bea{3.5}
 \frac{d}{dt}\io u\ln u +\io \frac{|\nabla u|^2}{u}
  &=&\chi \io \nabla u\cdot \nabla v +\xi \io\nabla u\cdot \nabla
  w\nn\\
  & & +\mu \io u(1+\ln u) (1-u-w) \qquad \mbox{for all $t\in (0,
  T)$}.
 \eea
Once more integrating by parts, in light of the Young inequality we
have
 \bea{3.6}
\chi \io \nabla u\cdot \nabla v
 &=& -\chi \io u \Delta v\nn\\
 &\le& \frac{1}{2}\io |\Delta v|^2 +\frac{\chi^2}{2}\io u^2 \qquad \mbox{for all $t\in (0,
  T)$}.
 \eea
Similarly,
 \bas
\xi \io\nabla u\cdot \nabla  w =-\xi \io u \Delta w=\xi \io u
(-\Delta w),
 \eas
which in view of Lemma \ref{lem_deltaw} entails that
 \bea{3.7}
 \xi \io \nabla u\cdot \nabla w
 &\le& \xi \|w_0\|_{L^\infty(\Omega)}\io uv +\xi K \io u\qquad \mbox{for all $t\in (0,
  T)$}
 \eea
because $u\ge 0$. Here we use the Young inequality and Lemma
\ref{lem3.2} to estimate the first term on the right
 \bas
\xi \|w\|_{L^\infty(\Omega)}\io uv
 &\le&\frac{\xi \|w_0\|_{L^\infty(\Omega)}}{2}\io u^2 +\frac{\xi \|w_0\|_{L^\infty(\Omega)}}{2}\io
 v^2\nn\\
 &\le& \frac{\xi \|w_0\|_{L^\infty(\Omega)}}{2}\io u^2 +c_1\qquad \mbox{for all $t\in (0,
  T)$}
 \eas
with $c_1:=\frac{\xi  \|w_0\|_{L^\infty(\Omega)}}{2}\cdot M(2)$,
whereas we employ Lemma \ref{mass} to deal with the second term on
the right of (\ref{3.7})
 \bas
\xi K \io u\le \xi K m^*     \qquad \mbox{for all $t\in (0,
  T)$}.
 \eas
Thus, we find that
 \be{3.8}
 \xi \io \nabla u\cdot \nabla w \le \frac{\xi \|w\|_{L^\infty(\Omega)}}{2}\io u^2 +c_2\qquad \mbox{for all $t\in (0,
  T)$},
 \ee
 where $c_2:=c_1+\xi K m^*$. As to the last term in (\ref{3.5}), by $u\ge 0, ~w\ge 0$, (\ref{bound}), Lemma \ref{mass} and the basic
 inequality $\max_{z\ge 0}(-z \ln z) =\frac{1}{e}$ we obtain
 \bea{3.9}
\mu \io u(1+\ln u) (1-u-w)
 &=&\mu \io u -\mu \io u^2 -\mu \io u w \nn \\
 & & +\mu \io u\ln u -\mu \io u^2
 \ln u +\mu \io (-u\ln u) w\nn \\
 &\le& \mu \io u +\mu \io u\ln u -\mu \io u^2
 \ln u +\mu \io (-u\ln u) w\nn \\
 &\le& \mu \io u\ln u -\mu \io u^2
 \ln u + \mu m^* +\frac{\mu}{e}\, \|w_0\|_{L^\infty(\Omega)} \cdot
 |\Omega|
 \eea
 for all $t\in (0, T)$. Collecting (\ref{3.6}), (\ref{3.8}) and
 (\ref{3.9}) along with (\ref{3.5}) leads to
 \be{3.10}
\frac{d}{dt}\io u\ln u +\io \frac{|\nabla u|^2}{u}\le \frac{1}{2}\io
|\Delta v|^2+\mu\io u\ln u +c_3 \io u^2-\mu\io u^2\ln u +c_4
 \ee
for all $t\in (0, T)$, where $c_3 :=\frac{\chi^2+\xi
\|w_0\|_{L^\infty(\Omega)}}{2}$ and $c_4 :=c_2 + \mu m^*
+\frac{\mu}{e}\, \|w_0\|_{L^\infty(\Omega)} \cdot
 |\Omega|$.\\
In order to cancel the first term on the right of (\ref{3.10}), we
test the second equation of (\ref{0}) by $-\Delta v$ and use the
Young inequality to find
 \bea{3.11}
\frac{1}{2}\frac{d}{dt} \io |\nabla v|^2 +\io |\nabla v|^2
 &=& -\io |\Delta v|^2 -\io u\Delta v\nn\\
 &\le& -\frac{1}{2}\io
|\Delta v|^2 +\frac{1}{2} \io u^2 \qquad \mbox{for all $t\in (0,
  T)$}.
 \eea
Adding this to (\ref{3.10}) yields
 \bas
\frac{d}{dt}\bigg\{\io u\ln u + \frac{1}{2} \io |\nabla v|^2\bigg\}
 + \io \frac{|\nabla u|^2}{u} + \io |\nabla v|^2
 \le \mu\io u\ln u +A \io u^2-\mu\io u^2\ln u +c_4
 \eas
for all $t\in (0, T)$, where $A: =\frac{1}{2}+c_3$. Adding $\io u\ln
u$ to both sides of this and dropping the nonnegative term $\io
\frac{|\nabla u|^2}{u}$ on the left, we find that $y(t):=\io u\ln u
+ \frac{1}{2} \io |\nabla v|^2$, $t\in (0, T)$, satisfies the
differential inequality
 \bas
y'(t) +y(t)\le \io \Big[(1+\mu) u\ln u +Au^2 -\mu u^2\ln u\Big]
+c_4,
 \eas
which in view of Lemma \ref{lem3.1} implies
 \bas
y'(t) +y(t)\le c_5 \qquad \mbox{for all $t\in (0,
  T)$}
 \eas
with $c_5:= L\cdot |\Omega|+c_4$. Upon ODE comparison, this yields
 \bas
y(t)\le \max \Big\{c_5, y(0)\Big\} \qquad \mbox{for all $t\in (0,
  T)$},
 \eas
which proves (\ref{3.3}) and (\ref{3.4}). \qed

\subsection{A bound for $\io u^2 +\io |\nabla v|^4$}
To build a bound for $\io u^2$,  we shall need the following
generalization of the Gagliardo-Nirenberg inequality for the general
case when $r>0$ (cf. \cite[Lemma A.5]{taowin8} for a detailed
proof), which extends the standard case when $r\ge 1$ in
\cite{biler_hebisch_nadzieja}.
\begin{lem}\label{generalized_GN}
Let $ \Omega\subset \R^2$ be a bounded domain with smooth boundary,
and let $p\in (1,\infty)$ and $r \in (0,p)$. Then there exists $C>0$
such that for each $\eta>0$ one can pick $C_\eta>0$ with the
property that
 \be{GN}
    \|u\|_{L^p(\Omega)}^p \le \eta \|\nabla u\|_{L^2(\Omega)}^{p-r}
        \big\|u \ln |u| \big\|_{L^r(\Omega)}^r
    + C\|u\|_{L^r(\Omega)}^p + C_\eta
\ee holds for all $u \in W^{1,2}(\Omega)$.
\end{lem}
By applying (\ref{2.11}) to $p=2$ and using Lemma \ref{lem3.2} we
establish an energy inequality involving $\io u^2$.

\begin{lem}\label{lem3.5}
Let $n=2$, $T\in (0, \tm)$, $\chi>0,~ \xi>0$ and $\mu>0$ , and
assume (\ref{init}). Then there exists $C_1>0$ independent of $T$
such that the solution of (\ref{0})-(\ref{02}) satisfies
\be{3.13}
      \frac{d}{dt} \io u^2 +\io |\nabla u|^2 \le  \chi^2 \io
      u^2|\nabla v|^2 -\mu \io u^3 +C_1
      \qquad \mbox{for all } t \in
      (0,T).
  \ee
\end{lem}
\proof We apply (\ref{2.11}) to $p=2$ to obtain
 \be{3.14}
      \frac{d}{dt} \io u^2 +\io |\nabla u|^2
     \le \chi^2  \io  u^2|\nabla v|^2 + 2\xi \|w_0\|_{L^\infty(\Omega)} \io u^2 v
       +2 (\mu+\xi K)\io  u^2 -2\mu \io u^3
       \ee
for any $p>1$ and each $t\in (0, T)$. Here we invoke the Young
inequality and Lemma \ref{lem3.2} to estimate
 \bea{3.15}
2\xi \|w_0\|_{L^\infty(\Omega)} \io u^2 v
 &\le& \frac{\mu}{2} \io u^3 + \frac{128}{27\mu^2}
 \Big(\xi\|w_0\|_{L^\infty(\Omega)}\Big)^3 \io v^3\nn\\
 &\le& \frac{\mu}{2} \io u^3 +c_1 \qquad\textrm{for all $t\in (0,
 T)$}
 \eea
with $c_1: =\frac{128}{27\mu^2}
(\xi\|w_0\|_{L^\infty(\Omega)})^3\cdot M(3)$, where $M(3)$ is
defined by Lemma \ref{lem3.2}. Similarly, we have
 \bas
2 (\mu+\xi K)\io  u^2 \le \frac{\mu}{2} \io u^3 +c_2
\qquad\textrm{for all $t\in (0,
 T)$}
 \eas
with $c_2: =\frac{128}{27\mu^2} (\mu+\xi K)^3\cdot |\Omega|$. This
in conjunction with (\ref{3.14}) and (\ref{3.15}) leads to
(\ref{3.13}) with $C_1:=c_1+c_2$. \qed
In order to deal with the first integral term on the right of
(\ref{3.13}), we further derive the following energy inequality for
$\io |\nabla v|^4$.
\begin{lem}\label{lem3.6}
Let $n\in \{2, 3\}$, $T\in (0, \tm)$, $\chi>0,~ \xi>0$ and $\mu>0$ ,
and assume (\ref{init}). Then there exists $C_2>0$ independent of
$T$ such that the solution of (\ref{0})-(\ref{02}) fulfills
  \bea{3.16}
    \hspace*{-8mm}
  &&  \frac{d}{dt} \io |\nabla v|^{4} +\io |\nabla v|^4
   + \io  \Big| \nabla |\nabla v|^2 \Big|^2
    \le 2 \int_{\pO} |\nabla v|^{2} \frac{\partial |\nabla v|^2 }{\partial \nu}
    +(n+4)\io u^2|\nabla  v|^{2}
    \eea
  for all  $t \in (0,T)$.
\end{lem}
\proof The proof is based on straightforward computations using the
second equation in (\ref{0}), and it was actually proved in
\cite[the proof of Lemma 3.3; see (3.12)-(3.13) therein]{taowin2}.
Thus, we prevent us from repeating the details here. \qed

\begin{cor}\label{cor3.7}
Let $n=2$, $T\in (0, \tm)$, $\chi>0,~ \xi>0$ and $\mu>0$, and assume
(\ref{init}). Then there exists $C_2>0$ independent of $T$ such that
the solution of (\ref{0})-(\ref{02}) carries the property
  \bea{3.17}
    \hspace*{-8mm}
   \frac{d}{dt}\bigg\{\io u^2 &+& \io |\nabla v|^{4}\bigg\}
   +\io u^2+ \io |\nabla v|^{4} +\io |\nabla u|^2
   + \io  \Big| \nabla |\nabla v|^2 \Big|^2    \nn\\
    &\le&  2 \int_{\pO} |\nabla v|^{2} \frac{\partial |\nabla v|^2 }{\partial
    \nu} +(\chi^2 +6) \io u^2|\nabla
   v|^{2}  +C_2 \qquad \textrm{  for all  $t \in (0,T)$.}
    \eea
\end{cor}

\proof Adding (\ref{3.16}) to (\ref{3.13}) yields
  \bas
    \hspace*{-8mm}
   \frac{d}{dt}\bigg\{\io u^2 &+& \io |\nabla v|^{4}\bigg\}
   +\io |\nabla v|^{4} +\io |\nabla u|^2 + \io  \Big| \nabla |\nabla v|^2 \Big|^2\nn\\
    &\le&  2 \int_{\pO} |\nabla v|^{2} \frac{\partial |\nabla v|^2 }{\partial
    \nu} + (\chi^2 +6)\io u^2|\nabla
    v|^{2} -\mu\io u^3 +C_1
     \eas
  for all
$t\in (0, T)$, where $C_1$ is provided by Lemma 3.5 . Adding $\io
u^2 $ to both sides of this and using the inequality
 \bas
\io u^2 \le \mu\io u^3+\frac{4}{27\mu^2}\cdot |\Omega|
 \eas
thanks to the Young inequality, we obtain (\ref{3.17}) with $C_2:=
 C_1+\frac{4}{27\mu^2}\cdot |\Omega|$. \qed
In the two-dimensional setting, we shall show that the two integrals
on the right of (\ref{3.17}) can be cancelled by $\io |\nabla u|^2 +
\io  | \nabla |\nabla v|^2|^2$ on the left, which thereby results in
a bound for $\io u^2+\io |\nabla v|^4$.
\begin{lem}\label{lem3.8}
Let $n=2$, $T\in (0, \tm)$, $\chi>0,~ \xi>0$ and $\mu>0$ , and
assume (\ref{init}). Then there exists $C>0$ independent of $T$ such
that the solution of (\ref{0})-(\ref{02}) enjoys the property
  \begin{eqnarray}
    \io u^2 &\le & C \qquad \textrm{ for all  $t \in (0,T)$} \quad
    \textrm{and}\label{3.18}\\
  \io |\nabla v|^4 &\le & C \qquad \textrm{ for all  $t \in (0,T)$}. \label{3.19}
    \eea
\end{lem}
\proof Starting from (\ref{3.17}), we first estimate $\int_{\pO}
|\nabla v|^{2} \Big|\frac{\partial |\nabla v|^2 }{\partial
\nu}\Big|$.
This boundary-related integral has been mainly studied in
\cite[(3.15)]{ishida_seki_yokota},  and accordingly we have
 \be{3.21}
2\int_{\pO} |\nabla v|^{2}  \bigg|\frac{\partial |\nabla v|^2
}{\partial \nu} \bigg|\le \frac{1}{2} \io \Big| \nabla |\nabla v|^2
\Big|^2 +c_1
 \qquad\mbox{for all $t\in (0, T)$}.
 \ee
We next deal with $\io u^2|\nabla v|^{2}$.
For any $\eta>0$, Young's inequality yields
 \be{3.23}
(\chi^2+6)\io u^2|\nabla  v|^2 \le \eta \io |\nabla v|^6
+\frac{(\chi^2+6)^{\frac{3}{2}}}{\sqrt{\eta}} \io u^3 \qquad
\mbox{for all } t \in
  (0,T).
 \ee
Here we use the Gagliardo-Nirenberg inequality and (\ref{3.4}) to
estimate
 \bea{3.24}
 \eta \io |\nabla v|^{6} &=& \eta \Big\| |\nabla
 v|^2\Big\|_{L^{3}(\Omega)}^{3}\nn \\
  &\le& \eta c_{2} \Big\| \nabla |\nabla
 v|^2\Big\|_{L^{2}(\Omega)}^{2} \cdot
 \Big\| |\nabla  v|^2\Big\|_{L^1(\Omega)} +c_{2}
 \Big\| |\nabla  v|^2\Big\|_{L^1(\Omega)}^3\nn\\
 &\le& \eta c_{3} \Big\| \nabla |\nabla
 v|^2\Big\|_{L^{2}(\Omega)}^2 +c_{3}
 \qquad \mbox{for all } t \in
  (0,T)
 \eea
and invoke Lemma \ref{generalized_GN} along with (\ref{3.3}) and
Lemma \ref{mass} to handle
 \bea{3.25}
\frac{(\chi^2+6)^{\frac{3}{2}}}{\sqrt{\eta}} \io u^3 &=&
\frac{(\chi^2+6)^{\frac{3}{2}}}{\sqrt{\eta}}
\|u\|_{L^3(\Omega)}^3\nn\\
 &\le& \frac{(\chi^2+6)^{\frac{3}{2}}}{\sqrt{\eta}}\Big[\eta \|\nabla u\|_{L^2(\Omega)}^2
 \cdot\|u\ln u\|_{L^1(\Omega)}
 +c_{4}\|u\|_{L^1(\Omega)}^3 +c_{5}(\eta)\Big] \nn\\
 &\le& c_{6}\sqrt{\eta} \|\nabla u\|_{L^2(\Omega)}^2 +c_{7}(\eta)
 \qquad \mbox{for all } t \in
  (0,T).
 \eea
Taking $\eta>0$ sufficiently small fulfilling
$\eta\le\min\{\frac{1}{2c_{3}}, \frac{1}{c_{6}^2}\}$, from
(\ref{3.23})-(\ref{3.25}) we infer that
 \be{3.26}
 (\chi^2+6)\io u^2|\nabla  z|^2 \le \io |\nabla u|^2 +\frac{1}{2} \io \Big| \nabla |\nabla
v|^2 \Big|^2 +c_{8}\qquad \mbox{for all } t \in
  (0,T).
 \ee
Thus, from (\ref{3.17}), (\ref{3.21}) and (\ref{3.26}) we obtain
that $y(t) :=\io u^2+ \io |\nabla v|^{4}$, $t\in (0, T)$, satisfies
the differential inequality
 \be{3.27}
y'(t) +y(t)\le c_{9}
 \ee
with $c_{9}:= c_1+c_8+C_2$. Upon an ODE comparison, this entails
 \bas
y(t)\le \max\Big\{y(0), c_{9}\Big\}, \qquad \textrm{for all $t\in
(0, T)$},
 \eas
which implies (\ref{3.18}) and (\ref{3.19}). \qed
Lemma \ref{lem3.8} results in the following useful corollary that
will be used in the proof of Lemma \ref{lem3.10} below.
\begin{cor}\label{cor3.9}
Let $n=2$, $T\in (0, \tm)$, $\chi>0,~ \xi>0$ and $\mu>0$ , and
assume (\ref{init}). Then there exists $C_3>0$ independent of $T$
such that the solution of (\ref{0})-(\ref{02}) possesses the
property
  \be{3.28}
      \| v (\cdot, t)\|_{L^\infty(\Omega)} \le C_3 \qquad \textrm{ for all  $t \in (0,T)$}.
    \ee
Moreover, for any $4\le q< \infty$ there exists $M_1(q)>0$ such that
the solution of (\ref{0}) satisfies
  \be{3.29}
      \io |\nabla v|^q \le M_1(q) \qquad \textrm{ for all  $t \in (0,T)$}.
    \ee
\end{cor}
\proof (\ref{3.19}) in conjunction with Lemma \ref{lem3.2} leads to
 \bas
 \| v(\cdot, t)\|_{W^{1, 4}(\Omega)} \le c_1 \qquad \textrm{ for all  $t \in
 (0,T)$}.
 \eas
This, along with the Sobolev embedding $W^{1,
4}(\Omega)\hookrightarrow C^0(\bar{\Omega})$ thanks to $4>n=2$,
yields (\ref{3.28}). As to (\ref{3.29}), it immediately follows from
(\ref{3.18}) and the standard parabolic regularity theory (cf.
\cite[Lemma 4.1]{horstmann_winkler} or \cite[Lemma 1]{kowszy}). \qed
\subsection{A bound of $u$ in $L^p(\Omega)$}
\begin{lem}\label{lem3.10}
Let $n=2$, $T\in (0, \tm)$, $\chi>0,~ \xi>0$ and $\mu>0$ , and
assume (\ref{init}). Then for any $p>2$ there exists $C(p)>0$
independent of $T$ such that the solution of (\ref{0})-(\ref{02})
fulfills
  \be{3.30}
  \io u^p \le C(p)  \qquad \textrm{for all $t \in (0,T)$}.
   \ee
\end{lem}
\proof Starting from (\ref{2.11}) once again and neglecting the
nonnegative term $\frac{p-1}{2} \io u^{p-2} |\nabla u|^2$ on the
left, we arrive at
  \bas
 \frac{d}{dt} \io u^p \le
      \frac{p(p-1) \chi^2}{2}
       \io  u^p|\nabla v|^2 +p\xi \|w_0\|_{L^\infty(\Omega)} \io u^p v
       +p(\mu+\xi K)\io  u^p -p\mu \io u^{p+1}
 \eas
for any $p>2$ and each $t\in (0, T)$. Adding $\io u^p$ to both sides
of this and using (\ref{3.28}) we see that
 \be{3.31}
 \frac{d}{dt} \io u^p +\io u^p\le
 \frac{p(p-1) \chi^2}{2}  \io  u^p|\nabla v|^2
       + c_1\io  u^p -p\mu \io u^{p+1}
 \ee
for any $p>2$ and each $t\in (0, T)$, where $c_1:=1+p(\xi
\|w_0\|_{L^\infty(\Omega)}\cdot C_3 +\mu +\xi K)$. Here we invoke
the Young inequality and (\ref{3.29}) to estimate
 \bea{3.32}
\frac{p(p-1) \chi^2}{2}  \io  u^p|\nabla v|^2
 &\le& \frac{p\mu}{2}\io u^{p+1} +c_2\io |\nabla v|^{2(p+1)}\nn\\
 &\le& \frac{p\mu}{2}\io u^{p+1} +c_2\cdot M_1(2p+2)
 \eea
with some $c_2>0$ and $M_1(\cdot)$ defined by Corollary
\ref{cor3.9}. Similarly, we have
 \be{3.33}
c_1\io u^p \le \frac{p\mu}{2}\io u^{p+1} +c_3
 \ee
with some $c_3>0$. Collecting (\ref{3.31})-(\ref{3.33}) yileds
$y(t): =\io u^p$, $t\in (0, T)$, satisfies the differential
inequality
 \bas
y(t)+y(t)\le c_4
 \eas
where $c_4:=c_2\cdot M_1(2p+2) +c_3$. Upon an ODE comparison, this
yields
 \bas
y(t)\le \max\Big\{y(0), c_4\Big\} \qquad \textrm{for all $t\in (0,
T)$},
 \eas
which leads to (3.30).\qed
\begin{cor}\label{cor3.11}
Let $n=2$, $T\in (0, \tm)$, $\chi>0,~ \xi>0$ and $\mu>0$ , and
assume (\ref{init}). Then there exists $C>0$ independent of $T$ such
that the solution of (\ref{0})-(\ref{02}) possesses the property
  \be{3.34}
  \|v(\cdot, t)\|_{W^{1,\infty}(\Omega)}\le C  \qquad \textrm{for all $t \in (0,T)$}.
   \ee
\end{cor}
\proof (\ref{3.34}) is a direct consequence of (\ref{3.30}) for a
fixed $p>2$ and the standard parabolic regularity theory (cf.
\cite[Lemma 4.1]{horstmann_winkler} or \cite[Lemma 1]{kowszy}). \qed

\mysection{Boundedness. Proof of Theorem \ref{theo1}}
Although (\ref{3.34}) shows that $\nabla v(\cdot, t)$ is bounded in
$L^\infty(\Omega)$, $\nabla w(\cdot, t)$ might become unbounded in
$L^\infty(\Omega)$ in light of (\ref{nablaw}). Therefore, we cannot
directly apply the result of the well-known Moser-Alikakos iteration
\cite{alikakos} to the first equation in (\ref{0}) to gain the
boundedness of $u(\cdot, t)$ in $L^\infty(\Omega)$. To bypass $w$,
our strategy is to use (\ref{2.11}) as a starting point for our
proof.
\begin{lem}\label{lem4.1}
Let $n=2$, $T\in (0, \tm)$, $\chi>0,~ \xi>0$ and $\mu>0$ , and
assume (\ref{init}). Then there exists $C>0$ independent of $T$ such
that the solution of (\ref{0})-(\ref{02}) satisfies
  \be{4.1}
  \|u(\cdot, t)\|_{L^\infty(\Omega)}\le C  \qquad \textrm{for all $t \in (0,T)$}.
   \ee
\end{lem}
\proof We begin with (\ref{2.11})
 \bas
      \frac{1}{p}\frac{d}{dt} \io u^p +\frac{p-1}{2}\io  u^{p-2}|\nabla u|^2
      & \le &\frac{(p-1) \chi^2}{2}
       \io  u^p|\nabla v|^2 +\xi \|w_0\|_{L^\infty(\Omega)} \io u^p v \nn\\
       & &+(\mu+\xi K)\io  u^p -\mu \io u^{p+1}
 \eas
for any $p>1$ and each $t\in (0, T)$. Adding $\io u^p$ to both sides
of this and invoking (\ref{3.34}), we obtain
 \be{4.2}
\frac{d}{dt} \io u^p +\io u^p +\io |\nabla u^{\frac{p}{2}}|^2 \le
c_1 p^2 \io u^p
 \ee
for any $p\ge 2$ and each $t\in (0, T)$, where $c_1>0$, as all
subsequently appearing constants $c_2, c_3, \cdots$ are independent
of $T$ as well as of $p\ge 2$. We now use the Gagliardo-Nirenberg
inequality to deal with the last integral
 \bas
\io u^p=\|u^{\frac{p}{2}}\|_{L^2(\Omega)}^2
 \le c_2 \|\nabla u^{\frac{p}{2}}\|_{L^2(\Omega)}\cdot \|u^{\frac{p}{2}}\|_{L^1(\Omega)}
  +\|u^{\frac{p}{2}}\|_{L^1(\Omega)}^2
 \eas
for any $p\ge 2$ and each $t\in (0, T)$. By Young's inequality, this
yields
 \bas
 c_1 p^2 \io u^p
 &\le& \io |\nabla u^{\frac{p}{2}}|^2  +c_3
 p^4\|u^{\frac{p}{2}}\|_{L^1(\Omega)}^2\\
 &=&  \io |\nabla u^{\frac{p}{2}}|^2  +c_3
 p^4\bigg(\io u^{\frac{p}{2}}\bigg)^2 \qquad \textrm{for all $t\in
 (0, T)$}.
 \eas
Hence, (\ref{4.2}) entails that
 \bas
\frac{d}{dt} \io u^p +\io u^p \le c_3
 p^4\bigg(\io u^{\frac{p}{2}}\bigg)^2 \qquad \textrm{for all $t\in
 (0, T)$}.
 \eas
Upon integration, this shows that
 \bas
\io u^p \le \io u_0^p +c_3 p^4 \int_0^t e^{-(t-\tau)} \bigg(\io
u^{\frac{p}{2}}(\cdot, \tau)\bigg)^2 d\tau.
 \eas
Writing $p_k :=2^k$ and
 \bas
B_k :=\max_{t\in (0, T)}\io u^{p_k} (\cdot, t)
 \eas
for $k\in\{1, 2, \cdots\}$, we see that
 \be{4.3}
B_k\le |\Omega|\cdot \|u_0\|_{L^\infty(\Omega)}^{p_k} +c_3
p_k^4B_{k-1}^2 \qquad \textrm{for all $k\ge 1$},
 \ee
where we have used the simple fact that $\int_0^t e^{-s} ds\le 1$.
Now if $p_k^4B_{k-1}^2\le  \|u_0\|_{L^\infty(\Omega)}^{p_k}$ for
infinitely many $k\ge 1$, we have
 \bas
\sup_{t\in (0, T)} \bigg(\io u^{p_{k-1}}(\cdot,
t)\bigg)^{\frac{1}{p_{k-1}}}\le \bigg(\frac{
\|u_0\|_{L^\infty(\Omega)}^{p_k}}{p_k^4}\bigg)^{\frac{1}{2p_{k-1}}}
= \frac{\|u_0\|_{L^\infty(\Omega)}}{p_k^{\frac{4}{p_k}}}
\longrightarrow \|u_0\|_{L^\infty(\Omega)}\qquad\mbox{as
$k\to\infty$},
 \eas
which implies that
 \bas
\sup_{t\in (0, T)} \|u(\cdot, t)\|_{L^\infty(\Omega)} \le
\|u_0\|_{L^\infty(\Omega)}
 \eas
and thereby proves the lemma in this case.\\
Conversely, if $p_k^4B_{k-1}^2>  \|u_0\|_{L^\infty(\Omega)}^{p_k}$
for all sufficiently large $k$, then (\ref{4.3}) yields some $c_4>0$
such that
 \bas
B_k\le c_4 p_k^4B_{k-1}^2 \qquad \mbox{for all $k\ge 1$}.
 \eas
In view of the definition of $p_k$, this implies
 \bas
B_k\le c_4 {(16)}^k B_{k-1}^2  \le a^k B_{k-1}^2 \qquad \mbox{for
all $k\ge 1$}
 \eas
with $a:=(\max\{c_4, 16\})^2$. Thus, by induction we obtain
 \be{4.4}
B_k\le a^{k+\sum_{j=1}^{k-1} 2^j (k-j)} \cdot B_0^{2^k} \qquad
\mbox{for all $k\ge 1$}.
 \ee
Here we observe that
 \bas
k+\sum_{j=1}^{k-1} 2^j (k-j) &=&2+2^2+\cdots +2^k-k \\
 &\le & 2^{k+1}
\qquad \mbox{for all $k\ge 1$}.
 \eas
From this and (\ref{4.4}) we infer
 \bas
 B_k^{\frac{1}{p_k}} \le a B_0\cdot a ^{\frac{1}{p_k}} \qquad \mbox{for all $k\ge
 1$},
 \eas
which after taking $k\to \infty$ implies that
 \bas
\sup_{t\in (0, T)}\|u(\cdot, t)\|_{L^\infty(\Omega)} \le aB_0
 \eas
and thereby yields the assertion in this case. \qed

We are now in a position to prove Theorem \ref{theo1}.\abs
\proofc of Theorem \ref{theo1}. \quad
The statement of global classical solvability and boundedness is a
straightforward consequence of Lemma \ref{lem_loc} and Lemma
\ref{lem4.1}. \qed
\vspace*{10mm}
{\bf Acknowledgment.}
  Y. Tao is supported by the National Natural Science Foundation of China (No.
  11171061) and by Innovation Program of Shanghai Municipal Education Commission (No. 13ZZ046).

\begin{thebibliography}{99}
%
\bibitem{alikakos}
  \sc Alikakos, N.D.: \it $L^p$ bounds of solutions of reaction-diffusion equations.
  \rm  Comm.~Partial Differential Equations {\bf 4}, 827-868 (1979)

\bibitem{biler_hebisch_nadzieja}
  \sc Biler, P., Hebisch, W., Nadzieja, T.: \it The Debye system: Existence and large time behavior of solutions.
  \rm Nonlinear Analysis, TMA {\bf 23} (9), 1189-1209 (1994)

\bibitem{chaplain_lolas_2005}
  \sc Chaplain, M.A.J., Lolas, G.:
  \it Mathematical modelling of cancer cell invasion of tissue: the
   role of the urokinase plasminogen activation system.
  \rm Math. Models Methods Appl. Sci. {\bf 18}, 1685每1734 (2005)

\bibitem{chaplain_lolas_2006}
  \sc Chaplain, M.A.J., Lolas, G.:
  \it Mathematical modelling of cancer invasion of tissue: dynamic
  heterogeneity.
  \rm Net. Hetero. Med. {\bf 1}, 399-439 (2006)

\bibitem{Corrias1}
  \sc  Corrias, L.,  Perthame, B., Zaag, H:
  \it A chemotaxis model motivated by angiogenesis.
  \rm  C. R. Acad. Sci. Paris, Ser. I. {\bf 336}, 141-146 (2003)

\bibitem{Corrias2}
  \sc  Corrias, L.,  Perthame, B., Zaag, H:
  \it Global solutions of some chemotaxis and angiogenesis systems in high
      space  dimensions.
  \rm  Milan J. Math. {\bf 72}, 1-28 (2004)

\bibitem{herrero_velazquez1997_pisa}
  \sc Herrero, M.~A., Vel\'azquez, J.~J.~L.:
  \it A blow-up mechanism for a chemotaxis model.
  \rm Ann.~Scuola Normale Superiore {\bf 24}, 633-683 (1997)

\bibitem{horstmann_winkler}
  \sc Horstmann, D., Winkler, M.:
  \it Boundedness vs. blow-up in a chemotaxis system.
  \rm J. Diff. Eqns. {\bf 215}, 52-107 (2005)

\bibitem{ishida_seki_yokota}
  \sc Ishida, S., Seki, K., Yokota T.:
  \it Boundedness in quasilinear Keller-Segel systems of parabolic-parabolic type on non-convex bounded
  domains.
  \rm J. Differential Equations {\bf 256}, 2993-3010 (2014)

\bibitem{kowszy}
  \sc Kowalczyk, R., Szyma\'{n}ska, Z.:
  \it On the global existence of solutions to an aggregation model.
  \rm J. Math. Anal. Appl. {\bf 343}, 379-398 (2008)

\bibitem{LM2}
  \sc Li\c{t}canu, G.,  Morales-Rodrigo, C.:
  \it Asymptotic behaviour of global solutions to a model of cell
  invasion.
  \rm Math. Mod. Meth. Appl. Sci. {\bf 20}, 1721-1758 (2010)

\bibitem{liutao}
  \sc Liu, D., Tao, Y.:
  \it Global boundedness in a fully parabolic attraction-repulsion
  chemotaxis model.
  \rm Math. Methods Appl. Sci., to appear

\bibitem{nagai2001}
  \sc Nagai, T.: \it Blowup of Nonradial Solutions to Parabolic-Elliptic Systems Modeling
  Chemotaxis in Two-Dimensional Domains.
  \rm J.~Inequal.~Appl. {\bf 6}, 37-55 (2001)

\bibitem{OTYM}
  \sc Osaki, K., Tsujikawa, T., Yagi, A., Mimura, M.:
  \it Exponential attractor for a chemotaxis-growth system of equations.
  \rm Nonlinear Anal. TMA {\bf 51}, 119-144 (2002)

\bibitem{tao_jmaa09}
  \sc Tao, Y.:
  \it Global existence of classical solutions to a combined
    chemotaxis每haptotaxis model with logistic source.
  \rm J. Math. Anal. Appl. {\bf 354},  60-69 (2009)

\bibitem{tao_wang_non}
 \sc Tao, Y., Wang, M.:
 \it Global solution for a chemotactic每haptotactic model of
   cancer invasion.
 \rm Nonlinearity {\bf 21}, 2221每2238 (2008)

\bibitem{taowin2}
  \sc Tao, Y., Winkler, M.: \it Boundedness in a quasilinear parabolic-parabolic Keller-Segel system with subcritical sensitivity.
  \rm J. Differential Equations {\bf 252}, 692-715 (2012)

\bibitem{taowin8}
  \sc Tao, Y., Winkler, M.: \it Energy-type estimates and global solvability
  in a two-dimensional chemotaxis-haptotaxis model with remodeling of non-diffusible attractant.
  \rm J. Differential Equations {\bf 257}, 784-815 (2014)

\bibitem{taowin_non}
 \sc Tao, Y., Winkler, M.:
 \it Dominance of chemotaxis in a chemotaxis每haptotaxis model.
 \rm Nonlinearity {\bf 27}, 1225-1239 (2014)

\bibitem{taowin8b}
  \sc Tao, Y., Winkler, M.: \it Boundedness and stabilization in a multi-dimensional chemotaxis-haptotaxis model.
  \rm Proceeding of the Royal Society of Edinburg, Section: A
  Mathematics, to appear

 \bibitem{WW}
 \sc  Walker, C., Webb, G.F.:
 \it Global existence of classical solutions for a haptotaxis model.
 \rm SIAM J. Math. Anal. {\bf 38}, 1694-1713 (2007)

\bibitem{win_cpde2010}
 \sc Winkler, M.:
 \it Boundedness in the higher-dimensional
  parabolic-parabolic chemotaxis system with logistic source.
 \rm  Commun. Partial Differential Equations {\bf 35}, 1516-1537
 (2010)

\bibitem{win_jde}
  \sc Winkler, M.:
  \it Aggregation vs. global diffusive behavior in the higher-dimensional Keller-Segel model.
  \rm J. Differential Equations {\bf 248}, 2889-2905 (2010)

\bibitem{win_arX2011}
  \sc Winkler, M:
  \it Finite-time blow-up in the higher-dimensional
  parabolic-parabolic Keller-Segel system.
  \rm Journal de Math\'ematiques Pures et Appliqu\'ees {\bf 100}, 748-767 (2013), {\tt arXiv:1112.4156v1}

\bibitem{win_jde2014}
  \sc Winkler, M.:
  \it Global asymptotic staility of constant equilibria in a fully
  parabolic chemotaxis system with strong logistic dampening.
  \rm J. Differential Equations {\bf 257}, 1056-1077 (2014)

%
\end{thebibliography}
\end{document}